\documentclass{gtart_h}

\def\ifplaintex{\expandafter\ifx\csname documentclass\endcsname\relax}

\def\gtp{{\mathsurround=0pt\it $\cal G\mskip-2mu$eometry \&\ 
$\cal T\!\!$opology $\cal P\!$ublications}}  

\def\Addressesr{\bigskip
{\small \parskip 0pt \leftskip 0pt \rightskip 0pt plus 1fil \def\\{\par}
\sl\theaddress\par
\medskip
\rm Email:\stdspace\tt\theemail\hfill\rm Received:\qua\receiveddate \par}}

\def\recd{{\small Received:\qua\receiveddate\ifx\reviseddate\relax
\else\qquad Revised:\qua\reviseddate\fi\par}} 

\def\lognumber#1{\def\thelognumber{#1}}
\def\volumenumber#1{\def\thevolumenumber{#1}}
\def\volumeyear#1{\def\thevolumeyear{#1}}
\def\papernumber#1{\def\thepapernumber{#1}}
\def\pagenumbers#1#2{\def\startpage{#1}\def\finishpage{#2}}
\def\published#1{\def\publishdate{#1}}

\def\received#1{\def\receiveddate{#1}}

\def\accepted#1{\def\accepteddate{#1}}

\long\def\asciiabstract#1{\long\def\theasciiabstract{#1}}

\let\\\par\let\thelognumber\relax\let\thevolumenumber\relax
\let\thepapernumber\relax\let\thevolumeyear\relax\let\startpage\relax
\let\finishpage\relax\let\publishdate\relax\let\receiveddate\relax
\let\reviseddate\relax\let\accepteddate\relax\let\theasciititle\relax
\let\theasciiauthors\relax
\let\theasciiabstract\relax

\let\theasciiemail\relax

\ifplaintex
\font\logobig=cmssbx10 scaled 3836
\font\logomed=cmssbx10 scaled 2557
\else
\font\logobig=cmssbx10 scaled 4200
\font\logomed=cmssbx10 scaled 2800
\fi

\long\def\makeagttitle{   
\count0=\startpage
\agt\hfill      
\hbox to 45truept{\vbox to 0pt{\vglue -13truept{\logomed A\kern -.37em{\logobig 
T}\kern -.38em G}\vss}\hss}
\break
{\small Volume \thevolumenumber\ (\thevolumeyear)
\startpage--\finishpage\nl
Published: \publishdate}

\vglue .25truein

{\parskip=0pt\leftskip 0pt plus
1fil\def\\{\par\smallskip}{\Large\bf\thetitle}\par\medskip} \vglue
0.05truein

{\parskip=0pt\leftskip 0pt plus 1fil\def\\{\par}{\sc\theauthors}
\par\medskip}%
 
\vglue 0.03truein 

{\small\leftskip 25truept\rightskip 25truept{\bf Abstract}\stdspace\theabstract

{\bf AMS Classification}\stdspace\theprimaryclass
\ifx\thesecondaryclass\relax\else; \thesecondaryclass\fi\par
{\bf Keywords}\stdspace \thekeywords\par}\vglue 7truept

}   

\ifplaintex
\font\phead=cmsl9 scaled 950
\font\pnum=cmbx10 scaled 913
\font\pfoot=cmsl9 scaled 950
\headline{\vbox to 0pt{\vskip -4.5mm\line{\small\phead\ifnum
\count0=\startpage ISSN 1472-2739 (on-line) 1472-2747 (printed)
\hfill {\pnum\folio}\else\ifodd\count0\def\\{ }%
\ifx\theshorttitle\relax\thetitle\else\theshorttitle\fi\hfill{\pnum\folio}
\else\def\\{ and }{\pnum\folio}\hfill\ifx\theshortauthors\relax\theauthors
\else\theshortauthors\fi\fi\fi}\vss}}
\footline{\vbox to 0pt{\vglue 0mm\line{\small\pfoot\ifnum\count0=\startpage
\copyright\ \gtp\hfill\else
\agt, Volume \thevolumenumber\ (\thevolumeyear)\hfill\fi}\vss}}
\else
\font=cmsl9 scaled 1050
\font\lnum=cmbx10 
\font=cmsl9 scaled 1050
\makeatletter
\def\@oddhead{{\small\lhead\ifnum\count0=\startpage ISSN 1472-2739 
(on-line) 1472-2747 (printed)\hfill {\lnum\number\count0}\else\ifodd\count0
\def\\{ }\ifx\theshorttitle\relax \thetitle \else\theshorttitle\fi\hfill
{\lnum\number\count0}\else\def\\{ and }{\lnum\number\count0}
\hfill\ifx\theshortauthors\relax 
\theauthors\else\theshortauthors\fi\fi\fi}}\def\@evenhead{\@oddhead}
\def\@oddfoot{\small\lfoot\ifnum\count0=\startpage\copyright\ \gtp\hfill\else
\agt, Volume \thevolumenumber\ (\thevolumeyear)\hfill\fi}
\def\@evenfoot{\@oddfoot}
\makeatother
\fi
\let\maketitlepage\makeagttitle
\let\makeshorttitle\maketitlepage
\let\maketitle\maketitlepage

\newwrite\gtoutfile
\long\gdef\makeheadfile{  
{\def\\{, }\def\s{ }
\immediate\openout\gtoutfile head.xxx
\immediate\write\gtoutfile{Proxy-for: \ifx\theasciiauthors\relax
\theauthors\else\theasciiauthors\fi\s<\ifx\theasciiemail\relax\theemail\else\theasciiemail\fi>}
\immediate\write\gtoutfile{\noexpand\\}
\immediate\write\gtoutfile{Authors: \ifx\theasciiauthors\relax
\theauthors\else\theasciiauthors\fi}
{\def\\{ }\immediate\write\gtoutfile{Title: \ifx\theasciititle\relax
\thetitle\else\theasciititle\fi}}
\immediate\write\gtoutfile{Subj-class: GT or SG, GR etc}
\immediate\write\gtoutfile{MSC-class: \theprimaryclass\ifx\thesecondaryclass\relax\else, \thesecondaryclass\fi}
\immediate\write\gtoutfile{Journal-ref: Algebr. Geom. Topol. \thevolumenumber\s
(\thevolumeyear) \startpage-\finishpage}
\immediate\write\gtoutfile{Comments: Published by Algebraic and
Geometric Topology at}
\immediate\write\gtoutfile{\s\s\s  http://www.maths.warwick.ac.uk/agt/AGTVol\thevolumenumber/agt-\thevolumenumber-\thepapernumber.abs.html}
\immediate\write\gtoutfile{\noexpand\\}
\immediate\write\gtoutfile{}
\ifx\theasciiabstract\relax
\immediate\write\gtoutfile{\theabstract}\else
\immediate\write\gtoutfile{\theasciiabstract}\fi
\immediate\write\gtoutfile{}
\immediate\write\gtoutfile{\noexpand\\}
\immediate\write\gtoutfile{}
\immediate\closeout\gtoutfile}}  

\def\maketitlepage{\makeagttitle\makeheadfile}
\let\makeshorttitle\maketitlepage
\let\maketitle\maketitlepage

\lognumber{32}
\volumenumber{5}
\volumeyear{2005}
\papernumber{32}
\pagenumbers{751}{768}
\received{16 June 2004} 
\accepted{24 June 2005}
\published{23 July 2005}

\usepackage{amsmath,amssymb}

\theoremstyle{plain}
\newtheorem*{Main}{Main Theorem}

\newtheorem{Thm}{Theorem}
\newtheorem{Cor}[Thm]{Corollary}
\newtheorem{Lem}[Thm]{Lemma}
\newtheorem{Prop}[Thm]{Proposition}

\newtheorem*{Cl}{Claim}

\theoremstyle{definition}
\newtheorem*{Def}{Definition}

\newcommand{\interior}{^{ \kern-5pt ^\circ}}
\newcommand {\bd}{\partial}

\newcommand {\R}{{\mathbb R}}

\newcommand {\HH}{{\mathbb H}}
\newcommand {\nb}{\text{Nbh}}
\newcommand {\st}{\text{Stab}}

\newcommand {\ld}{\Lambda}

\newcommand {\cA}{ {\cal A}}
\newcommand {\cB}{ {\cal B}}

\newcommand {\cC}{{\cal C}}
\newcommand {\cE}{{\cal E}}
\newcommand {\cF}{{\cal F}}

\newcommand {\cP}{{\cal P}}

\newcommand {\cY}{{\cal Y}}

\newcommand {\diam}{{\rm diam}}

\begin{document}
\title {Bootstrapping in convergence groups}
\author{Eric L. Swenson}

\address{Mathematics Department, Brigham Young University\\Provo, 
 UT 84604, USA }

\email{eric@math.byu.edu}

\begin{abstract} We prove a true bootstrapping result for
convergence groups acting on a Peano continuum.

We give an example of a Kleinian group $H$ which is the amalgamation
of two closed hyperbolic surface groups along a simple closed curve.
The limit set $\Lambda H$ is the closure of a ``tree of circles"
(adjacent circles meeting in pairs of points).  We alter the action of
$H$ on its limit set such that $H$ no longer acts as a convergence
group, but the stabilizers of the circles remain unchanged, as does
the action of a circle stabilizer on said circle.  This is done by
first separating the circles and then gluing them together backwards.
\end{abstract}

\asciiabstract{%
We prove a true bootstrapping result for
convergence groups acting on a Peano continuum.
We give an example of a Kleinian group H which is the amalgamation
of two closed hyperbolic surface groups along a simple closed curve.
The limit set Lambda H is the closure of a `tree of circles'
(adjacent circles meeting in pairs of points).  We alter the action of
H on its limit set such that H no longer acts as a convergence
group, but the stabilizers of the circles remain unchanged, as does
the action of a circle stabilizer on said circle.  This is done by
first separating the circles and then gluing them together backwards.}


\maketitle

\section{Introduction}
When a group $G$ acts by homeomorphisms on a compact metric space
$X$, then, in general, little can be said about $G$ or $X$.
However if $G$ acts on $X$ as a convergence group then we can often
say much more.  For instance if $X$ is a topological circle then $G$
``is" a Fuchsian group where the action of $G$ on $X$ is given by the
Fuchsian action on the circle at infinity of $G$ \cite{TUK}, \cite{GAB},
\cite{CAS-JUN}. For any metric space, if $G$ acts on $X$ as a uniform
convergence group, then $G$ is word hyperbolic and $X= \bd G$
\cite{BOW}.  More generally if $G$ acts on $X$ as a geometrically
finite convergence group, then $G$ is relatively hyperbolic with
boundary $X$ \cite{YAM}.

The purpose of this note is to examine what conditions will allow
us to prove that $G$ acts on $X$ as a convergence group. In
particular we are interested in ``bootstrapping", that is using
subgroups (with presumably smaller limit sets) which act as
convergence groups to prove that $G$ acts as a convergence group
on $X$.   We will restrict ourselves to that case where $X$ is a
Peano continuum without cut points.  This restriction  is
justified to some extent by the results of \cite{BES-MES},
\cite{BOW1}, \cite{BOW2} and \cite{BOW3}. This is a continuation of
a program began in \cite{SWE1}, \cite{SWE2}.  Our goal (which we
achieve) is a natural bootstrapping result when $X$ is an
$n$-sphere ($n>1$) and the subgroups have limit sets
$(n-1)$-spheres in $X$.  This follows from our main theorem.
(Recall that a map is {\em essential} if it is not homotopic to a
constant map)

\begin{Main}
Let $X$ be a Peano continuum, without cut points, which does not
admit an essential map to the circle. If $G$ acts on $X$ by
homeomorphisms and $\cA$ is a $G$-invariant collection of
connected closed subsets of $X$, then $G$ acts as a convergence
group on $X$ provided the following conditions are satisfied:
\begin{itemize}
\item $ \cA$ is null, that is for any $\epsilon >0$,
the set of elements of $\cA$ with diameter at least $\epsilon$,
$\{A \in \cA: \diam(A) > \epsilon \}$, is finite.
\item $ \cA$ is fine, that is for any $x,y \in X$
there exists a finite $\cB \subset \cA$ such that $\cup \cB$
separates $x$ from $y$.
\item $\st (A)$ acts as a convergence group on $A$ for each $A \in
\cA$.
\end{itemize}
\end{Main}

We also have results which apply when the elements of $\cA$ are
not connected, but they are more technical.

 \section{Background}
\begin{Def} Let $G$ be a group acting by homeomorphisms on a space
$X$, a sequence $(g_i)\subset G$ acts as a {\em convergence
sequence on $X$} if $\exists\, n,p \in X$ such that for any compact
$C \not \ni n$, $g_i(C) \to p$. We call  $n$ the {\em repeller}
and $p$ the {\em attractor} of $(g_i)$ in its action on $X$.
 (It is possible that $n=p$).
 It is easily shown that if $(g_i)$ acts as a
convergence sequence on $X$ with attractor $p$ and repeller $n$,
then $(g_i^{-1})$ acts as a convergence sequence on $X$ with
attractor $n$ and repeller $p$. We say that $G$ acts as a
convergence group on $X$ if every sequence of distinct elements of
$G$ has a subsequence acting as a convergence sequence on $X$.
\end{Def}
A {\em continuum} is a compact connected metric space, and a {\em
Peano continuum} is a locally connected continuum.  For $A \subset
X$ and $B,C$ connected subsets of $X-A$ we say $A$
separates $B$ from $C$ if $B$ and $C$ lie in different components
of $X-A$. We say that $A$ separates $Y \subset X$ if $A$ separates
two points of $Y$.

\begin{Def}
Let $A$ and $B$ be closed sets of a Peano Continuum $X$.  We say
$A$ and $B$ {\em cross\/} if either $A\cap B  \neq \emptyset$
or $A$ separates $B$ and $B$ separates $A$.

A  sequence $A_1, \dots A_n$ is called a {\em crossing sequence}
 from $A_1$ to $A_n$ if $A_i$ crosses $A_{i+1}$ for all $1\le i<n$.
The sequence is called a minimal crossing sequence if $A_i$
doesn't cross $A_j$ for $|i-j|>1$.

\end{Def}

\begin{Def} Let $X$ be a Peano continuum, and
 $G$ be a group which acts by homeomorphisms on $X$.  If  $\cA$ is a
$G$-invariant collection of closed subsets of $X$, we say that the
pair $(G, \cA)$ is a fine pairing on $X$ if the following
conditions are satisfied:
\begin{enumerate}
\item $\cA$ is { \em cross connected}, that is for any
$A,B \in \cA$ there is a crossing sequence in $\cA$ from $A$ to
$B$.
\item $\cA$ is {\em null}.  That is: For any $\epsilon >0$,
the set of elements of $\cA$ with diameter at least $\epsilon$,
$\{A \in \cA: \diam(A) > \epsilon \}$ is finite.
\item $\cA$ is {\em fine}.  That is: For any $x,y \in X$
there exists a finite $\cB \subset \cA$ such that $\cup \cB$
separates $x$ from $y$.
\end{enumerate}
\end{Def}
The following is the main result of \cite{SWE2}.
\begin{Thm}\label{T:old}
 Let $X$ be a Peano continuum without cut
 points, and $(G, \cA)$ be a fine pairing on $X$.  If, for each
$A \in \cA$, $\st (A)$ acts as a convergence group on $X$, then
  $G$ acts as a convergence group on $X$.
  \end{Thm}

  The difficulty with applying this result is that in order to
  prove that $G$ acts as a convergence group on $X$ we must prove
  that each element of an infinite family of subgroups acts as a
  convergence group
  on $X$.  If, in the above result, we could replace the phrase
   ``$\st (A)$ acts as a convergence group on $X$"
   with the phrase ``$\st (A)$ acts as a convergence group on $A$",
   then we would be able to bootstrap.
 This was what the author set out to prove.
   The following counterexample explains why he failed.

\section{Counterexample}
The idea of this counterexample is quite simple, but the details
are not.  What we will do is start with a Kleinian group $H$ whose
limit set is a ``tree" of circles in $S^2$ (where a pair of
circles which intersect one another do so in a pair of points
which are the endpoints of a hyperbolic element). In fact $H$ will
be the free product of two surface groups amalgamated along an
element corresponding to a simple closed curve in both.  We now
take the tree of circles, pull out every other circle in it and
glue them back in interchanging each pair of  points in an
intersection. The action of $H$ on $\ld H$ will  result in an
action on the new tree of circles, but the action of $H$ on this
new tree of circles will not be a convergence action even though
it satisfies all of the conditions of the conjecture that the
author wished to prove.

This example is constructed using only what the author considers
to be the most elementary of techniques at his disposal.  An
expert in Kleinian groups could complete this construction much
faster, and without any need for the geometric group theory
arguments  which the author employs.

Let  $F$ be the Kleinian group generated by reflection in the
faces of a regular right-angled hyperbolic dodecahedron $D$. Let
$P$ (purple) and $Y$ (yellow) be hyperbolic planes which contain
adjacent faces (that is to say faces which intersect in an edge)
of $D$. Let $F_P = \st (P)$, the stabilizer of $P$ in $F$, and
similarly define $F_Y= \st (Y)$. Clearly $P/ F_P$ is a (purple)
closed hyperbolic orbifold and similarly $Y/F_Y$ is a (yellow)
closed hyperbolic orbifold. The intersection line $P \cap Y$
descends to a closed geodesic in both $P/F_P$ and $Y/ F_Y$, in
fact both $P/F_P$ and $Y/ F_Y$ are regular right angled hyperbolic
pentagon reflection orbifolds, and the intersection line will
descend to an edge of the orbifold.
 Notice that there is an isomorphism $\chi:P/F_P  \to Y/ F_Y$
 of these orbifolds which
 sends the intersection edge in $P/F_P$ to the intersection edge in
$Y/ F_Y$.

  Fix $d >0$. Using classical
 Fuchsian group results we can find a finite index
subgroup $H_P <F_P$, so that letting  $H_Y = \chi(H_P) < F_Y$ we
have:
\begin{enumerate}
\item The quotients $P/H_P$ and $Y/H_Y$ are orientable closed
surfaces. \label{s}
\item The intersection line $P \cap Y$ descends to the simple closed
 intersection curve in each of $P/H_P$ and $Y/H_Y$. \label{c}
\item $H_P \cap H_Y = \langle g \rangle$ where $\langle g \rangle$
is the stabilizer of $P\cap Y$ in both $H_P$ and $H_Y$.\label{nbs}
\item Every closed curve on $P/H_P$ or $Y/H_Y$ which intersects the
 intersection curve transversely has length at least $d$.
 \label{d}
\end{enumerate}
 Let  $H =\langle H_P,H_Y
 \rangle$, the subgroup of $F$ generated by $H_P$, and $H_Y$.

 We will show that for $d \gg 0$, $H = H_P \ast_{\langle g\rangle} H_Y$,
 the free product of $H_P$ with $H_Y$ amalgamated over $\langle g\rangle$.
   Let $\cP = \{h(P) |\, h \in H
 \}$ be the purple planes (translates of $P$ by $H$) and $\cY =\{
 h(Y)|\, h \in H \} $  be the yellow planes (translates of $Y$).

 We show that $P$ is not a translate of $Y$ in $H$ by showing
 that $P$ is not a translate of $Y$ in our original reflection
  group $F$.  To do this give each face of $D$ its own color, and
  tile $\HH^3$ with copies of $D$ by reflecting in the faces of
  $D$.  Notice that after reflection in a single face of $D$, the
  resulting 17-hedron (2 of the pentagon faces will be in the
  interior of the solid, and 10 pentagons will be amalgamated into
  5 hexagon faces) will have each face of one color.  It
  follows that each plane in the tiling of $\HH^3$ will have only
  one color (each plane will be tiled with pentagons of the same
  color).  Thus $P$ is not a translate of $Y$ in $F$ or in $H$, so
 $\cP \cap \cY = \emptyset$.

 We construct a bipartite graph $T$ with alternating purple and
 yellow vertices in the obvious way.  Let $T^{(0)} = \cP \cup \cY$
 so the vertices are the purple and yellow planes, and we connect two
 vertices with an edge if and only if the corresponding planes
 intersect.  The group $H$ acts on $T$ with out fixed points.
\begin{Def}
Let $x_0,x_1, \dots x_m \in \HH^3$. The union of the geodesic
segments $ \alpha = [x_0,x_1] \cup [x_1,x_2] \cup [x_2, x_3] \dots
\cup [x_{n-1}, x_n]$ is called a {\em broken geodesic} from $x_0$
to $x_n$, with break points $x_1, \dots x_{n-1}$.  The length
$\ell(\alpha)=\sum \ell([x_{i-1},x_i])$ is the sum of the lengths
of the segments. The angle between the segments $[x_{i-1}, x_i]$
and $[x_i,x_{i+1}]$ is called the break angle of $\alpha$ at
$x_i$.
\end{Def}
In \cite{CAN1}, Cannon proved the following theorem.
\begin{Thm} For any $\theta>0$ there exists $d(\theta) >0$ and
$\epsilon(\theta) >0$ with the following property: For any broken
geodesic $\alpha = \bigcup\limits_{i=1}^n[x_{i-1},x_i]$, with
break angles at least $\theta$,  satisfying $\ell([x_{i-1},x_i])
\ge d(\theta)$ for all $1<i<n$ and for any $a,b \in \alpha$, then
$d(a,b) \ge \epsilon(\theta) \ell(\alpha')$ where $\alpha'$ is the
subpath of $\alpha$ from $a$ to $b$.
\end{Thm}
When $d \ge d(\frac \pi 2)$, using Cannon's result we show that
$T$ is a tree. Suppose not. There is a minimal circuit in $T$, and
from this we construct a broken geodesic closed curve $\alpha
\subset \cup(\cP \cup \cY) \subset \HH^3$ where each geodesic
segment in $\alpha$ is $\alpha \cap B$ where $B$ is one of the
vertex planes of the circuit.  Homotopying $\alpha$ relative to
the intersection lines of our vertex planes, we may assume that
$\alpha$ has minimal length, that is to say no homotopy of
$\alpha$ relative to the intersection lines has shorter length.
(This argument is made very easy by the fact that the intersection
lines diverge from one another exponentially).

  Notice that every
segment of $\alpha$ has length at least $d$ by \eqref{d}.  We show
that no break angle in $\alpha$ is acute. Consider a break point
$q \in \alpha$, lying on the intersection of planes $P'\in \cP$,
$Y'\in \cY$. Since $P' \perp Y'$, there is a unique hyperbolic
plane $Q\ni q$ with $Q \perp P'$ and $Q\perp Y'$.

If the interiors of the 2 segments of $\alpha$ which meet at $q$
lie on same side of $Q$ (aside from the endpoint $q$) then we can
shorten $\alpha$ by sliding $q$ to that side of $Q$, and if not
then the break angle is not acute.
  Thus the
break angle of $\alpha$ at $q$ is not acute.

 Every segment is at
least $d$ long, so by Cannon's theorem $d(q,q) \ge \epsilon(\frac
\pi 2) \ell(\alpha)
>0$. This is a contradiction, so $T$ is a tree.
 A ping pong argument now shows that the stabilizer of $P$ in $H$ is $H_P$ and
 similarly the stabilizer of $Y$ in $H$ is $H_Y$.
 Thus $H = H_P \ast_{\langle g\rangle} H_Y$

Consider the Poncaire ball model for $\HH^3$, the ball of radius 1
about $\mathbf 0$.  Position $D$ so that ${\mathbf 0} \in P \cap Y$.
We move our perspective to the boundary $S^2= \{ {\mathbf x} \in \R
^3|\, ||{\mathbf x}|| =1 \}$ of the ball model for $\HH^3$. Each
hyperbolic plane in $\HH^3$ corresponds to a circle in $S^2$. Thus
the purple planes correspond to purple circles in $S^2$ and
similarly the yellow planes correspond to yellow circles. We will
denote by $\pi$ the boundary circle of $P$, and by $\upsilon$ the
boundary circle of $Y$. Two planes intersect in a line if and only
if the corresponding circles intersect in two points. Let $\Pi$ be
the collection of all the purple circles ($\Pi = \{ \bd \breve P
|\, \breve P \in \cP \}$), and similarly let $\Upsilon$ be the
collection of yellow circles ($\Upsilon = \{ \bd \breve Y |\,
\breve Y \in \cY \}$).

We now redefine the vertex set of $T$ to be $\Pi \cup \Upsilon$
and two vertices are connected by an edge if the vertex circles
intersect in two points.  Of course $T$ is still (the same)
bipartite (yellow and purple) tree.
$$X= \cup (\Pi \cup \Upsilon) \subsetneq \Lambda H \subset
S^2.\leqno{\rm Let}$$ 
As we will see, $X$ is a finite diameter path metric space.
For this we need the following result of negative curvature due to
Gromov \cite{GRO}, Cannon \cite{CAN1}, \cite{CAN2} and
others.
\begin{Prop} \label{P:quasi} For any $ \epsilon > 0$ there is a $\delta >0$ such that
if $p:[a,b] \to \HH^n$ is a path satisfying $ d(p(x),p(y))  \ge
\epsilon\ell(p([x,y])) $then $p$ and the geodesic from $p(a)$ to
$p(b)$ fellow travel, $[p(a),p(b)] \subset \nb(p([a,b]), \delta)$
and $ p([a,b]) \subset \nb([p(a),p(b)], \delta)$.
\end{Prop}
By increasing $d$, we can assume that $d \gg \delta$ (increasing
$d$ doesn't decrease $\epsilon$) where $\delta$ is the constant of
Proposition \ref{P:quasi}.

 Let $\alpha$ be a geodesic in $T$ ,
and define $\cC(\alpha)$ to be the collection of vertex circles of
$\alpha$.  For any circle $C \subset S^2 \subset \R^3$ the
circumference of  $C$ will be denoted $circ(C)$.

\begin{Lem} \label{L:geo}There is $S>0$ such that for any geodesic
$\alpha $ in $T$,$$\sum\limits_{C \in \cC(\alpha)}circ(C) \le S[
circ(B)]$$ where $B$ is the vertex of $\alpha$ closest to the
vertices $\{\pi, \upsilon \}$ in $T$.
\end{Lem}
\begin{proof}If suffices to prove the following:
There is a $\phi <1$, such that for any $B$, a vertex circle of
$T$, and $C$, an adjacent vertex circle of $T$ with $B$ closer (in
$T$) to $\{ \pi, \upsilon\}$ than $C$, then $circ(C) \le circ(B)
\phi $.

Using high school geometry, we can show that the Euclidean radius
of a circle $E$ on $S^2$ is $e^{-h}$ where $h$ is the hyperbolic
distance from $\mathbf 0$ to the hyperbolic plane (Euclidian
hemisphere in this model) whose boundary is $E$.

Let $S_C$ and $S_B$ be hyperbolic planes in $\HH^3$ corresponding
to $C$ and $B$ respectively. Let $\gamma$ be the shortest geodesic
from $\mathbf 0$ to $S_C$.  Let $\beta$ be the shortest path in
$\cup (\cP \cup \cY)$ with the same endpoints as $\gamma$. By
Proposition \ref{P:quasi}, $\beta$ and $\gamma$ $\delta$-fellow
travel. Notice that $\ell(\beta \cap S_B) >d $

  It follows that $d({\mathbf 0}, S_B) + d- 3 \delta \ge d({\mathbf 0},
S_C)$.  Thus since $\frac d 2 < d- 3\delta$, letting $\phi =
e^{-\frac d 2}$ suffices. We have $circ(C) \le  circ (B)\phi$
since this is true for the radii.
\end{proof}

Recall that is $X$ is the ``tree" of the purple and yellow circles
in the limit set of the Kleinian group $H$, $ \ld(H) \subset S^2$.
By Lemma \ref{L:geo}, $X$ is a finite diameter geodesic metric
space on which $H$ acts by homeomorphisms.

By orienting the intersection curve on $P/ H_P$, we obtain a well
defined orientation on the intersection line $P \cap Y$.  Let
$\pi$ and $\upsilon$ be the purple and yellow circles
corresponding to $ P$ and $ Y$ respectively. We now denote $\pi
\cap \upsilon = \{+,-\}$ where $+$ is the point at infinity of the
line $P\cap Y$ in orientation direction, and $-$ is the other
endpoint of $P \cap Y$. Since $H_P$ is a torsion free hyperbolic
surface group, it follows that $+$ is not a translate of $-$ under
the action of $H_P$, and similarly under the action of  $H_Y$.
However the stabilizer of $P$ in $H$ is $H_P$, and the stabilizer
of $Y$ is $H_Y$, and so $+$ is not a translate of $-$ in $H$.

 For any $\breve \pi \in \Pi$ and $\breve
\upsilon \in \Upsilon$ which intersect there is $h \in H$ with
$\breve \pi = h(\pi)$ and $\breve \upsilon = h(\upsilon)$.  Thus
$\breve \pi \cap \breve \upsilon = \{ h(+), h(-) \}$. In short the
intersection of any two circles in $X$ consists of a positive
point $h(+)$ and a negative point $h(-)$.

Now consider the disjoint union $W = (\cup \Pi) \sqcup (\cup
\Upsilon)$, that is $W$ is the disjoint union of the yellow and
purple circles, so each path component of $W$ is a single circle
where $W$ is given the path metric. Now for $\breve \pi \in \Pi$
and $\breve \upsilon \in \Upsilon$ where $\breve \pi$ and $\breve
\upsilon$ intersect in $X$, define $+_{\breve \pi}( \breve
\upsilon) \in \breve \pi \subset W$ to be the positive
intersection point in $\breve\pi$, and $-_{\breve \pi}( \breve
\upsilon) \in \breve \pi \subset W$ to be the negative
intersection point in $\breve \pi$. Define $+_{\breve \upsilon}(
\breve \pi), \,-_{\breve \upsilon}( \breve \pi) \in\breve \upsilon
\subset W$ similarly.

Notice that $H$ acts by homeomorphism on $W$ and that $X$ is the
quotient space of $W$ by the equivalence relation $+_{\breve
\upsilon}( \breve \pi)= +_{\breve \pi}( \breve \upsilon)$ and
$-_{\breve \upsilon}( \breve \pi) = -_{\breve \pi}( \breve
\upsilon)$. Also the action of $H$ on $X$ is the quotient action
of the action of $H$ on $W$.

We now form a different quotient of $W$ on which $H$ acts. Define
the bassackward  equivalence relation $\sim$ by $+_{\breve
\upsilon}( \breve \pi)\sim  -_{\breve \pi}( \breve \upsilon)$ and
$-_{\breve \upsilon}( \breve \pi)\sim +_{\breve \pi}( \breve
\upsilon)$. Notice that since $+$ is not a translate of $-$, for
any $h \in H$ $h$ preserves the equivalence classes of $W$.  It
follows that $H$ acts via homeomorphism on the quotient space $Z=
W /\sim$. Where we give $Z$ the path metric calculating the arc
length of a path in each circle.

 Notice that the tree structure of $Z$ is exactly
the same as the tree structure of $X$  (the same pairs of purple
and yellow circles intersect), and so it follows from Lemma
\ref{L:geo} that $Z$ is a finite diameter path metric space.

Now for any $h \in H$, consider $h:S^2 \to S^2$.  Since $h$ is
conformal, then $h$ has continuous first derivative.  Since $S^2$
is compact we have that $h:S^2 \to S^2$ is Lipschitz with constant
$c>0$.  It follows that for any circle $C \subset S^2$, $h:C \to
h(C)$ is Lipschitz with constant $c$,{ \em even when $C$ and
$h(C)$ are given the path metric inherited from $S^2$}. Thus $h:Z
\to Z$ is also Lipschitz with constant $c$.  It follows that $h$
extends to a homeomorphism $h: \hat Z \to \hat Z$ where $\hat Z$
is the metric completion of $Z$.

\begin{Thm} $\hat Z$ is compact.
\end{Thm}
\begin{proof}It suffices to show that $Z$ is totally bounded (for every
$\epsilon >0$ there is a finite covering of $Z$ by
$\epsilon$-balls). This follows easily from the compactness of a
circle together with Lemma \ref{L:geo}, which implies that there
are only a finite number of circles above a given size. Cover the
finite collection of ``large" circles with a finite collection of
balls of very small radius. The collection of balls of radius
$\epsilon$ with the same centers will cover $Z$.
\end{proof}

Now the group $H$ acts by homeomorphism on the Peano continua
without cut-points $\hat Z$. Let $\cA = \Pi \cup \Upsilon$  (that
is the elements of $\cA$ are the circles).  The collection of
circles  $\cA$ is cross-connected, null and fine since this is
true when $\cA$ is consider as a collection of closed subset of
$\ld(H)$, and there is a homeomorphism of $\hat Z$ to $\ld (H)$
which permutes the elements of $\cA$.  Thus $(H, \cA)$ is a fine
pairing on $\hat Z$.  Notice that for any circle $A \in \cA$, the
stabilizer $\st(A)$ is the same as it was in $\ld H$ (a coFuchsian
subgroup) and the action of $\st(A)$ on $A$ is also the same (the
action of a coFuchsian group on it limit circle).  Thus for each
$A \in \cA$, $\st(A)$ acts on $A$ as convergence group. Now $(H,
\cA)$ on $\hat Z$ satisfies all the hypothesis of the failed
conjecture, but
\begin{Thm} The action of  $H$ on $\hat Z$ is not a convergence
action.
\end{Thm}
\begin{proof} Recall that  $\pi$ and $\upsilon$ are the circles
corresponding to the planes $P$ and $Y$ resp.  The stabilizer
$\st(\pi) = H_P$ and $\st(\upsilon) = H_Y$.  Recall that $H_P \cap
H_Y = <g>$ where $g$ represents the intersection curve of $P/H_P$
and $Y/H_Y$.  The sequence $(g^i)$ ($i>0$) acts as a convergence
sequence on both $\pi$ and $\upsilon$.  In its action on $\pi$,
the attractor of $(g^i)$ is the point $+_\pi(\upsilon) \in \pi$,
and the repeller is the point $-_\pi(\upsilon) \in \pi$.  On the
other hand, in its action on $\upsilon $, $(g^i)$ has attractor
$+_\upsilon(\pi) \in \upsilon$ and repeller $-_\upsilon(\pi) \in
\upsilon$.  However in $\hat Z$, $+_\pi(\upsilon) =-_\upsilon(\pi)
$ and $+_\upsilon(\pi) = -_\pi(\upsilon)$.

  Let $n \in \hat Z$, and choose $a \in \pi -\{n,-_\pi(\upsilon) \}$ and $b
\in \upsilon - \{n,-_\upsilon(\pi) \}$.  We have that $g^i(a) \to
+_\pi(\upsilon)$ and $g^i(b) \to -_\upsilon(\pi)$, and the same
would be true for subsequence of $(g^i)$.  Thus $n$ is not the
repeller of any convergence subsequence of $(g^i)$ (since there is
not a single attractor), and so $(g^i)$ has no convergence
subsequence and $H$ doesn't act as a convergence group on $\hat
Z$.
\end{proof}

\section{Bootstrapping}
\begin{Def} Let  $X$ be a Peano continuum without cut points, and
 $G$ be a group which acts by homeomorphisms on $X$.  Let $\cA$ be a
$G$-invariant collection of  closed subsets of $X$ such that the
following conditions are satisfied:
\begin{enumerate}
\item $|A| > 2$ for all $A\in \cA$.
\item $\cA$ is {\em null}.  That is: For any $\epsilon >0$,
the set of elements of $\cA$ with diameter at least $\epsilon$,
$\{A \in \cA: \diam(A) > \epsilon \}$ is finite.\label{null}
\item $\cA$ is {\em fine}.  That is: For any $x,y \in X$
there exists a finite $\cB \subset \cA$ such that $\cup \cB$
separates $x$ from $y$.\label{fine}
\item For each $A \in \cA$, the stabilizer of $A$, $\st(A) =\{g \in G : g(A)=A\}$
acts as a convergence group on $A$. \label{stab}
\item For any
$A,B \in \cA$  with $|A\cap B| \le 2$:
 \begin{itemize}
 \item If $A \cap B =\{c\}$, then for any $b \in B-\{c\}$
 there exists finite crossing sequence
$ A, A_1, \dots A_n, B \in \cA$ with $D \cap A_i = \emptyset$ for
all $1\le i \le n$, where $D = \{b,c\}$.
 \item If $|A\cap B| \neq 1$ then there is a  crossing sequence
$ A, A_1, \dots A_n, B \in \cA$ with $D \cap A_i = \emptyset$ for
all $1\le i \le n$, where $D= A\cap B$.
\end{itemize}
   \label{connect}
\end{enumerate}
If these conditions are satisfied then we say that the pair $(G,
\cA)$ is a {\em local convergence pairing} on $X$.  Clearly every
local convergence pairing is  a fine pairing.
\end{Def}
We will show that if $(G, \cA)$  is a local convergence pairing on
$X$, then $G$ acts as a convergence group on $X$. We translate the
following two results from \cite{SWE2}, the second of which is not
due to the author, into
 the setting we are interested in.

\begin{Lem} \label{L:pigion}  Let $(G, \cA)$ be a local
convergence pairing on $X$ and $\cB \subset \cA$ with $\cB$ finite
and cross-connected. If $(g_i)$ is a sequence of elements of $G$
then either $(g_i(\cup \cB))$ is null or there exists $B \in \cB$,
$A \in \cA$ and a subsequence $(g_{i_n})$ of $(g_i)$ with
$g_{i_n}(B) =A$ for all $n$.
\end{Lem}

\begin{Lem} \label{L:topo} If $M$ is a compact connected locally
connected Hausdorff space and a closed set $A$ doesn't separate
the nonempty connected set $B \subset M-A$ from the nonempty
connected set $C\subset M-A$, then there exists a neighborhood $V$
of $A$ which doesn't separate $B$ from $C$.
\end{Lem}

\begin{Lem}\label{L:small}
Let $Y$ be a Peano continuum without cut-points. For any point $x
\in Y$ and for any $\epsilon >0$ there is a $\delta > 0$ such that
$Y- B(x,\epsilon)$ is contained in a single component of $Y- B(x,
\delta)$.
\end{Lem}
\begin{proof} Suppose not, then there exists  $\epsilon >0$, such that
for any $\delta>0$, $Y- B(x, \epsilon) $ is not contained in a
single component of $Y- B(x, \delta)$.  In particular, for each $n
> \frac 1 \epsilon$, there exist $a_n, b_n \in Y- B(x, \epsilon) $ such that
$a_n$ and $b_n$ lie in different components of $Y- B(x, \frac 1
n)$.  By compactness, we may assume that $a_n \to a \in Y- B(x,
\epsilon)$ and $b_n \to b \in Y- B(x, \epsilon)$. Let $A$ and $B$
be  closed connect neighborhoods of $a$ and $b$ respectively with
$x \not \in A \cup B$.    There exists a $\delta>0$ such that
$A\cap B(x, \delta) = \emptyset=B\cap B(x, \delta)$ and (using
Lemma \ref{L:topo}) $B(x , \delta)$ does separate $A$ from $B$.
This contradicts the choice of $x_n$ and $y_n$ for $n\gg 0$.
\end{proof}

\begin{Lem} \label{L:good} Let $(G, \cA)$ be a local convergence pairing on
$X$, and suppose $B, A_1, \dots A_m \in \cA$ is a crossing sequence. Let
$(f_i)\subset G$ act as a convergence sequence on $B$ with
attractor $p$ and repeller $n$.  If for all $1 \le j \le m$, we
have that  $n, f_i^{-1}(p) \not \in A_j$ for all $i$, then
$f_i(A_j) \to p$ for all $1\le j \le m$.
\end{Lem}
\begin{proof} By nullity, since $p \not \in f_i(A_j)$, it suffices to show that
$d(p,f_i(A_j)) \to 0$.

 Proceeding by induction consider the case where $j=1$.
  If $B \cap A_1 \neq \emptyset$, then there is $a \in (B \cap
  A_1)- \{n\}$ and $f_i(a) \to p$, yielding $d(p, f_i(A_1)) \to 0$.
  Now consider the case where $B$ separates $A_1$ and $A_1$
  separates $B$.  Let $B^+$ be the intersection of $B$ with that
  component of $X-A_1$ which contains $p$, and let $B^-= B- B^+$.
 We know from \cite{TUK2} that $n$ is not an isolated point of
 $B$, so $B^\pm \neq \{n\}$. Let $b \in B^--\{n\}$ and $\hat b \in B^+
-\{n\}$.
 We have that $f_i(b), f_i(\hat b) \to
 p$, but $f_i(A_1)$ separates $f_i(b)$ from $f_i(\hat b)$.  It
follows from local connectivity of $X$ that $d(p, f_i(A_1)) \to
0$.

Given that $f_i(A_j) \to p$, one of the following is true:
\begin{itemize}
\item $\exists\, a \in A_j \cap A_{j+1}$ and so $f_i(a) \to p$.
\item $f_i(A_j)$ separates $f_i(A_{j+1})$ and by Lemma
\ref{L:small}, $d(p, f_i(A_{j+1})) \to 0$.
\end{itemize}
Either way $d(p, f_i(A_{j+1})) \to 0$ and so $f_i(A_{j+1}) \to p$.
\end{proof}

\begin{Lem} \label{L:nn} Let $(G, \cA)$ be a local convergence pairing on
$X$, and $A,B \in \cA$.  If $(f_i)\subset \st(A) \cap \st(B)$ acts
as a convergence sequence on both $A$ and $B$ with $p$ the
attractor of $(f_i)$ in its action on $A$ and on $B$ then the
repeller $n$ of $(f_i)$ in its action on $A$ is also the repeller
of $(f_i)$ in its action of $B$.
\end{Lem}
\begin{proof}
If $A\cap B \neq \{p\}$, then there is $a \in A\cap B -\{p \}$ and
the convergence sequence $(f_i^{-1})$ has attractor $n$ in its
action on $A$ and so $f_i^{-1}(a) \to n$. Thus the convergence
sequence $(f_i^{-1})$ has attractor $n$ in its action on $B$.  It
follows that the repeller of $(f_i)$ in its action on $B$ is $n$.

Now we deal with the case where $A\cap B = \{p\}$. Suppose that
the repeller of $(f_i)$ in its action on $B$ is $\hat n \neq n$.
Notice that $f_i(p) =p$ for all $i$.

 By
\eqref{connect} there is a finite crossing sequence $A,A_1, \dots
A_m, B$ with $p \not \in A_j$ for all $j$, and also passing to a
subsequence of $(f_i)$, we may assume that $f_i( n) \not \in A_j$
for all $i$ and $1\le j \le m$.

We now apply Lemma \ref{L:good} to $(f_i^{-1})$ ( which has
attractor $n$ and repeller $p$ in its action on $A$) and $A, A_1,
\dots A_m$ we find that $f_i^{-1}(A_m) \to n$.  Applying Lemma
\ref{L:good} to $(f_i^{-1})$ and the sequence $B, A_m$ we get that
$f_i^{-1}(A_m) \to \hat n$.  Thus $n= \hat n$.
\end{proof}

\begin{Thm} \label{T:match} Let $(G, \cA)$ be a local convergence pairing on
$X$, and $A,B \in \cA$ with $A \cap B \neq \emptyset$.  If
$(f_i)\subset G$ acts as a convergence sequence on both $A$ and
$B$ and $n$ is the repeller and $p$ the attractor of $(f_i)$ in
its action on $A$, then the attractor of $(f_i)$ in its action on
$B$ is $p$ and the repeller of $(f_i)$ in its action on $B$ is
$n$.
\end{Thm}
\begin{proof}
Suppose we have distinct $a,b \in A\cap B -\{n\}$,  with $f_i(a)
\to p$ and $f_i(b) \to p$.  It follows that $p$ is the attractor
of $(f_i)$ in its action on $B$, and so by Lemma \ref{L:nn} $n$ is
the repeller of $(f_i)$ in its action on $B$.

Similarly if we have distinct $a,b \in A\cap B -\{p\}$, then
$(f_i^{-1})$ has attractor $n$ and repeller $p$ in its action on
$B$, and so $p$ is the attractor, and $n$ the repeller of $(f_i)$
in its action on $B$.

We are left with the case where $A\cap B \subset \{n, p \}$.
Replacing $(f_i)$ with $(f_i^{-1})$ if need be (and so exchanging
the rolls of $n$ and $p$), we may assume that $n \in A \cap B$.

\rk{Case I\quad  $ p \in A\cap B$}

Notice that $(f_i)$ leaves $A\cap B$ invariant. By \eqref{connect}
there is a crossing sequence $A,A_1, \dots A_m,B$ with $p,n \not
\in A_j$ for all $j$. Since $f_i^{-1}(p) \not \in A_j$ for all $i$
and $j$, Lemma \ref{L:good} implies that $f_i(A_m) \to p$.

If $A_m \cap B= \emptyset$, then by Lemma \ref{L:good} applied to
$(f_i)$ and the crossing sequence $B, A_m$, $f_i(A_m)$ must
converge to the attractor of $(f_i)$ in its action on $B$, but
$f_i(A_m) \to p$, so $p$ is the attractor of $(f_i)$ in its action
on $B$ and Lemma \ref{L:nn} implies that $n$ is the repeller of
$(f_i)$ in its action on $B$.

 If not there exists $b \in A_m \cap B$. It follows that
  $f_i(b) \to p$. Since $(f_i)$ leaves that set $A \cap B$
invariant it follows from \cite{TUK} that the attractor and
repeller of $(f_i)$ in its action on $B$ are contained in the set
$A \cap B$. All together this implies that $p$ is the attractor of
$(f_i)$ in its action on $B$ and by Lemma \ref{L:nn}  that $n$ is
the repeller of $(f_i)$ in its action on $B$.

\rk{Case II\quad $p \not \in A\cap B$}

By \eqref{connect} there is a crossing sequence $A,A_1, \dots A_m,
B$ with $n \not \in A_j$ for all $j$, and also passing to a
subsequence of $(f_i)$, we may assume that $f_i^{-1}( p) \not \in
A_j$ for all $i$ and $j$. Thus we may apply Lemma \ref{L:good} to
show that $f_i(A_m) \to p$.

If $A_m \cap B \neq \emptyset$, then $f_i(A_m \cap B) \to p$,
so $p \in B$ which puts us back in Case I.

If $A_m \cap B = \emptyset$, then by Lemma \ref{L:good} applied to
$(f_i)$ and the crossing sequence $B, A_m$, we find that
$f_i(A_m)$ converges to the attractor of $(f_i)$ in its action on
$B$. It follows that $p$ is the attractor of $(f_i)$ in its action
on $B$, so $p \in B$ which puts us back in Case I.
\end{proof}
\begin{Lem} \label{L:close} Let $(G, \cA)$ be a local convergence pairing on $X$ and $A,B \in
\cA$. If $(g_i) \subset \st (A)$ acts as a convergence sequence on
$A$ with attractor $p$ and repeller $n$, and $A \cap B \neq
\emptyset$, then there exists a subsequence $(h_i) \subset (g_i)$
such that for any compact $C \subset B-\{n\}$, $h_i(C) \to p$.
\end{Lem}
\begin{proof} The proof has two cases.

\rk{Case I\quad $\diam (g_i(B)) \not \to 0$}

By Lemma \ref{L:pigion}
 there is $D \in A$ and a subsequence $(k_i)\subset (g_i)$
  with $k_i(B)=D$ for
all $i$.  Let $f_i =k_1^{-1}k_i \in \st(A) \cap \st(B)$.  Clearly
$(f_i)$ acts as a convergence sequence on $A$ with attractor
$k_1^{-1}(p)$ and repeller $n$.  Using \eqref{stab} we find a
subsequence $(q_i) \subset (f_i)$ with $(q_i)$ acting as a
convergence subsequence on $B$.  By Theorem \ref{T:match} $(q_i)$
has attractor $k_1^{-1}(p)$ and repeller $n$ in its action on $B$.
Let $h_i = k_1q_i$ so $(h_i) \subset (k_i)\subset (g_i)$. For any
compact $C \subset B-\{n\}$, by definition $q_i(C) \to
k_1^{-1}(p)$.  It follows that $h_i(C) \to p$.

\rk{Case II\quad $\diam (g_i(B)) \to 0$} 

If there is $a \in A \cap B -\{n\}$ then $g_i(a) \to p$ and so $g_i(B)
\to p$ and we are done.  Thus we need only consider the case where $A
\cap B = \{n \}$.  Using \eqref{connect} and passing to a subsequence
$(h_i) \subset (g_i)$ we obtain a crossing sequence $A, A_1, \dots
A_m,B$ with $n \not \in A_j$ for all $j$, and $h_i^{-1}( p) \not \in
A_j$ for all $i$ and $j$. Thus by Lemma \ref{L:good}, $h_i(A_m) \to
p$.

 If
$b \in A_m \cap B$ we get that $h_i(b) \to p$ and so $h_i(B) \to
p$.

If $A_m \cap B = \emptyset$, then for each $i$, $h_i(B)$ is
separated by the set $h_i(A_m)$.  It follows from Lemma
\ref{L:small} that $d(p, h_i(B)) \to 0$ and so $h_i(B) \to p$.
\end{proof}
\begin{Cor} \label{C:far} Let $(G, \cA)$ be a local convergence
pairing on $X$ and $A,B \in \cA$. If $(g_i) \subset \st (A)$ acts
as a convergence sequence on $A$ with attractor $p$ and repeller
$n$, then there exists a subsequence $(h_i) \subset (g_i)$ such
that for any compact $C \subset B-\{n\}$, $h_i(C) \to p$.
\end{Cor}
\begin{proof}
Using \eqref{connect} we find a minimal crossing sequence $A,A_1
\dots A_{m-1}, A_m=B \in  \cA$.
 If $A_1 \cap A = \emptyset$ then by Lemma
\ref{L:good}, $g_i(A_1) \to p$.
If $A_1 \cap A \neq \emptyset$ we apply Lemma \ref{L:close} and
pass to a subsequence of $(g_i)$ so that for any compact $C\subset
A_1 -\{n\}$, $g_i(C) \to p$.

Now by induction we may assume that for any compact $C \subset A_k
- \{n \}$, $g_i(C) \to p$.
If $a \in A_{k+1} \cap A_k$, then since $A_{k+1} \cap A =
\emptyset$, $g_i(a) \to p$ and it follows that $g_i(A_{k+1}) \to
p$.

Now consider when  $A_{k+1} \cap A_k = \emptyset$. Either
$g_i(A_k) \to p$ and by Lemma \ref{L:small} $g_i(A_{k+1}) \to p$.

 Otherwise
$\diam\, g_i(A_k) \not \to 0$, and arguing as in Case I of Lemma
\ref{L:close} we obtain a convergence sequence acting on $A_k$,
and applying Lemma \ref{L:good} to this convergence sequence and
the crossing sequence $A_k, A_{k+1}$ gives us a subsequence $(h_i)
\subset (g_i)$ with $h_i(A_{k+1}) \to p$.
\end{proof}

 \begin{Thm} \label{T:Main} If $(G,\cA)$ is a local convergence pairing on $X$ then
 $G$ acts as a convergence group on $X$.
 \end{Thm}
 \begin{proof}
 By Theorem \ref{T:old} it suffices to show that for each $A\in \cA$, $\st(A)$ acts as a
 convergence group on $X$.

 Suppose not, then there is $(f_i) \subset \st(A)$ with $(f_i)$ acting as a convergence sequence on
 $A$ with attractor $p$ and repeller $n$, but no subsequence of $(f_i)$ acts as convergence sequence
 on $X$ with attractor $p$ and repeller $n$.  Using the compactness of $X$, and passing to a
 subsequence, we may assume there is a sequence $(x_i) \subset X- \{n \}$ with $x_i \to x \neq n$
 and $f_i(x_i)\to \hat x \neq p$. Choose a finite $\cB \subset \cA$ such that $\cup \cB$ separates
 $n$ from $x$.  Using \eqref{connect} and Corollary \ref{C:far},
 passing to subsequences we may assume that $f_i(B) \to p$ for all $B \in \cB$,
  so $f_i(\cup \cB) \to p$.

  Let $U$ be the component of $X - \cup \cB$ which contains $n$.  Thus $\bd U \subset \cup \cB$ and
  so $f_i(\bd U) \to p$.  Notice that $f_i(A-U) \to p$ so $\diam(f_i(A-U)) \to 0$ and this implies that
  $\diam f_i(U \cap A)$ is bounded away from $0$ (since $|A| > 2$), and so $\diam(g_i(U))$ is bounded
  away from $0$.  Since $p$ is not a cut point one can show using Lemma \ref{L:topo} that
  $f_i(X-U) \to p$, but for all $i \gg 0$ $x_i \in X-U$. Thus $f_i(x_i)  \to p$ contradicting
  $\hat x \neq p$.
  \end{proof}

\begin{Lem}\label{L:H_1} Let $Z$ be a Peano continuum  which admits
no essential map to the circle, and let $A,B$ be disjoint closed
subsets of $Z$.  Let $c,d \in Z- [A\cup B]$. If $c$ and $d$ are
not separated by $A$ nor separated by $B$, then $c$ and $d$ are
not separated by $A \cup B$.
\end{Lem}
 Janiszewski proved this in the case where $Z$ was the
  2-sphere \cite[XVII 5.2]{DUG}.
\begin{proof}
Suppose not, that is $c$ and $d$ lie in different components of
$X-[A \cup B]$. Let $S^1 =I \cup J$ where $I$ and $J$ are arcs
with $I \cap J = \{\hat c, \hat d\}$.  Let $ \hat a \in
\,$Int$\,I$ and $\hat b \in \,$Int$\,J$.  Let $\hat U$ be the arc
of $S^1$ from $\hat a$ to $\hat b$ passing through $\hat c$, and
similarly defined $\hat V$ to be the arc of $S^1$ from $\hat a$ to
$\hat b$ passing through $\hat d$.

 The space $X$ is locally
path connected since it is a Peano continuum and so we may choose
a mapping $\alpha : I \to X$ with $\alpha(\hat c) =c$, and
$\alpha(\hat d) =d$ with the property that $\alpha(I) \cap B =
\emptyset$ (since $B$ doesn't separated). Similarly we choose a
map  $\beta:J \to X$ with $\beta (\hat c) =c$, and $\beta(\hat d)
=d$ with the property that $\beta(J) \cap A = \emptyset$. We paste
$\alpha$ and $\beta$ to obtain $\gamma: S^1 \to X$.

Let $X-[A\cup B] =U\cup V$ where $U$ and $V$ are open sets of $X$,
$U\cap V = \emptyset $ and $c \in U$ and $d \in V$. By the Tietze
Extension theorem, there exists a map $f:U \cup A \cup B \to \hat
U $ with $f(A) = \{\hat a\}$ and $f(B) = \{\hat b\}$.  Similarly
there exists a map $g: V \cup A \cup B \to \hat V$ with $g(A) =
\{\hat a\}$ and $g(B) = \{\hat b\}$.  We paste $f$ and $g$ to get
a continuous $h:X \to S^1$.  By construction $h \circ \gamma: S^1
\to S^1$ is homotopic to the identity map and so essential.
Clearly $h$ is essential since $h \circ \gamma$ is essential. Thus
$X$ admits an essential map, $h$, to $S^1$ which is a
contradiction.
\end{proof}
\begin{Cor} \label{C:compcon} Let $Z$ be a Peano continuum  which admits
no essential map to the circle, let $A$ be a closed subset of $Z$,
and let $c,d \in Z- A$. If $c$ and $d$ are  separated by $A$, then
$c$ and $d$ are separated by some component of $A$.
\end{Cor}
\begin{proof}  Let $(A_i)$ be a nested sequence ($A_{i+1} \subseteq A_i$)
of closed subsets of $A$ with the property that for each $i$,
$A_i$ separates $c$ from $d$.  By Lemma \ref{L:topo},
$\bigcap\limits_i A_i$ separates $c$ from $d$. Thus by Zorn's
Lemma, there exists a minimal closed subset  $B$ of $A$ which
separates $c$ from $d$. Since no proper closed subset of $B$
separates $c$ from $d$, by Lemma \ref{L:H_1}, $B$ is connected.
\end{proof}
\begin{Main}
Let $X$ be a Peano continuum, without cut points, which does not
admit an essential map to the circle. If $G$ acts on $X$ by
homeomorphisms and $\cA$ is a $G$-invariant collection of
connected closed subsets of $X$, then $G$ acts as a convergence
group on $X$ provided the following conditions are satisfied:
\begin{itemize}
\item $ \cA$ is null, that is for any $\epsilon >0$,
the set of elements of $\cA$ with diameter at least $\epsilon$,
$\{A \in \cA: \diam(A) > \epsilon \}$, is finite.
\item $ \cA$ is fine, that is for any $x,y \in X$
there exists a finite $\cB \subset \cA$ such that $\cup \cB$
separates $x$ from $y$.
\item $\st (A)$ acts as a convergence group on $A$ for each $A \in
\cA$.
\end{itemize}
\end{Main}
\begin{proof} First we remove from $\cA$ all singleton sets
(sets with only one element).  We must verify that this new $\cA$
satisfies the three original hypothesis.  Two of them are obvious,
and the fact that  $\cA$ is still fine follows from $X$ having no
cut points and Corollary \ref{C:compcon}.  Thus we may assume that
$\cA$ has no singletons, and since each element of $\cA$ is
connected, each element of $\cA$ is uncountable.

  We
must show that $(G,\cA)$ is a local convergence pairing on $X$. We
are given all of the conditions for a local convergence pairing
except \eqref{connect} which we state again here.
 For any $A,B \in
\cA$ with $|A\cap B| \le 2$:
 \begin{itemize}
 \item If $A \cap B =\{c\}$ we choose any $b \in B-\{c\}$ and let $D=
 \{b,c\}$.
 \item If $|A\cap B| \neq 1$ we let $D= A\cap B$.
 \end{itemize}
 there exists finite crossing sequence
$ A_1, A_2, \dots A_n \in \cA$ with $A_1 \cap A \neq \emptyset$,
and $A_n \cap B \neq \emptyset$  but $D \cap A_i = \emptyset$ for
all $1\le i \le n$.

There are three cases.
\rk{Case I\quad $A \cap B = \emptyset$}  

We must find a crossing sequence in $\cA$ from $A$ to $B$.

We start by showing that $\cup \cA$ is connected.  Suppose not,
then there are open sets $U,V$ of $X$ with $\cup \cA
\subset U \cup V$, $\cup \cA \cap U \cap V  = \emptyset$, $\cup
\cA \cap U \neq \emptyset \neq \cup \cA \cap V$.  By Corollary
\ref{C:compcon}, we may assume that $\bd U$ is connected.   Since $X$
has no cut points, there exist distinct $x,y \in \bd U$. By fine,
there exist a finite $\cF \subset \cA$ such that $\cup \cF$
separates $x$ and $y$.  Since $\cup \cF$ separates $x$ from $y$,
and $\bd U$ is connected, then $\cup \cF \cap \bd U \neq
\emptyset$.  This is a contradiction since $\cup \cF \subset U
\cup V$.  Thus $\cup \cA$ is connected.

We now show that $\cup \cA$ is cross connected.  We define the
cross-component of $E \in \cA$ as
$$\gimel(E) = \cup \{ C \in \cA:\text{ there is a crossing
sequence in }\cA \text{ from }E\text{ to }C\}$$ If $e \in E$, we
define $\gimel(e) = \gimel(E)$ (note that this is well defined).

\begin{Cl} For each $e \in \cup \cA$, there exists a
basis $\{U_\alpha\}$ of connected open sets at $e$ in $X$ with the
property that $\bd U_\alpha \subset \gimel(e)$ for all $\alpha$.
\end{Cl}

 Let $z\in X-\{e\}$ and
$y \in \gimel(e)-\{e,z\}$.  By fineness, there is a finite $\cF
\subset \cA$ so that the component $U$ of $X- \cup \cF$ which
contains $e$ has the property that $y,z \not \in \overline{U}$.
Since $e$ is not a cut point, using compactness and Lemma
\ref{L:topo}, we may assume that $ \overline{U}$ doesn't separate
$y$ from $z$. Thus $\bd U$ separates $e$ from $\{y,z\}$ which is
contained in a single component of $X- \bd U$.  By Corollary
\ref{C:compcon}, there is a single component $P$ of $\bd U$ which
separates $e$ from $y$ (of course $P$ also separates $e$ from
$z$). Since $e,y \in \gimel(e)$, which is connected, $P \cap
\gimel(e) \neq \emptyset$.  It follows that $P \subset \gimel(e)$
since $P$ is a connected subset of a finite union of elements of
$\cA$. Thus we have show that for every $z \in X-\{e\}$ there is a
connected open set $V\ni e$ with $\bd V \subset \gimel(e)$ but $z
\not \in \overline{V}$.  Using compactness, the Claim follows.

Let $a\in A$ and $b \in B$, and choose $U_a, U_b$, disjoint
connected open neighborhoods of $a$ and $b$ respectively, with
$\bd U_a \subset \gimel(a)$ and $\bd U_b \subset \gimel(b)$. Now
for every other $e \in \cup \cA$ choose
 $U_e$, connect open neighborhood of $e$, with the property that $\bd
U_e \subset \gimel(e)$ and $a,b \not \in U_e$. The
collection$\{U_e: e \in \cup \cA\}$ is an open cover of $\cup \cA$
Since  $\cup \cA$ is connected there exists a finite chain of open
sets  $ U_{e_1}, \dots U_{e_n},$ with $U_{e_i} \cap U_{e_{i+1}}
 \neq \emptyset$ for each $1\le i \le n$ where $a=
e_1$ and $b=e_n$. Notice that since $U_a \cap U_b = \emptyset$,
$n>2$. We may assume this chain is minimal.

 We now show that for any $1< i <n$, $\gimel(e_{i-1})= \gimel(e_i) =
 \gimel(e_{i+1})$. Since the chain is minimal, $U_{e_{i-1}} \cap
 U_{e_{i+1}} = \emptyset$,  $U_{e_{i-1}} -\overline{U_{e_{i}}} \neq \emptyset$
 and $U_{e_{i+1}} -\overline{U_{e_{i}}} \neq \emptyset$.  Thus
 $U_{e_{i-1}}-\bd U_{e_{i}}$ is not connected and
  and  neither is $U_{e_{i-1}}-\bd U_{e_{i}}$.  Since  $U_{e_{i-1}}$ and
   $U_{e_{i+1}}$ are
 connected sets, it follows that $\bd U_{e_i} \cap U_{e_{i-1}}
 \neq \emptyset \neq \bd U_{e_i} \cap U_{e_{i+1}}$.  Notice that $\bd
 U_{e_i} \subset \gimel(e_i)$ which is connected.  Since $U_{e_{i-1}} \cap
 U_{e_{i+1}} = \emptyset$, $\gimel(e_i) \cap \bd U_{e_{i-1}} \neq
 \emptyset \neq \gimel(e_i) \cap \bd U_{e_{i+1}}$  Thus
 $\gimel(e_i) \cap \gimel(e_{i-1}) \neq \emptyset \neq \gimel(e_i) \cap
 \gimel(e_{i+1})$.  It follows that $\gimel(e_{i-1})= \gimel(e_i) =
 \gimel(e_{i+1})$.

 Thus by induction, $\gimel(A)=  \gimel(a) =\gimel(b)= \gimel(B)$ and so
 there is a finite crossing sequence from $A$ to $B$ as required.

\rk{Case II\quad $|A\cap B| = \{c\}$}  

Choose some $b \in B- \{c \}$ and
some $a \in A - \{c\}$. Using compactness of $X$ and fineness it
can be shown that for any neighborhood $W$ of $c$, there exists a
neighborhood $U \subset W$ of $c$ and a finite subset $\cF \subset
\cA$ such that $\bd U \subset \cF$.   Using Lemma \ref{L:topo} can
find an open neighborhood $U$ of $c$  with $a,b $ in the same
component $N$ of $X- U$ and with $\bd U \subset \cup \cF$ where
$\cF$ is a finite subset of $\cA$.  Using the nullness of $\cA$, we
may assume that $a,b \not \in \cup \cF$. By Lemma \ref{L:H_1},
there is a component $\cup \cE$ of $\cup \cF$ which separates $c$
from $N$.   Thus $\cE$ is cross-connected and $\{b,c\} \cap \cup
\cE = \emptyset$.  Since $\cup \cE$ separates $a$ from $c$  and
separates $b$ from $c$, there is a crossing sequence in $\cE$ from
$A$ to $B$ missing the set $\{b,c\}$.

\rk{Case III\quad $|A \cap B| =2$} 

This is similar to (but easier than) the
previous case and will be left to the reader (simply find a finite
subset of $\cA$ which separates the two points of intersection).
\end{proof}

\Addressesr

\end{document}